\documentclass{amsart}

\usepackage{latexsym,amsfonts,amssymb,exscale,enumerate,comment}
\usepackage{amsmath,amsthm,amscd}

\usepackage{url}
\usepackage[bookmarks=true,%
    colorlinks=true,%
    linkcolor=blue,%
    citecolor=blue,%
    filecolor=blue,%
    menucolor=blue,%
    urlcolor=blue,%
    breaklinks=true]{hyperref}

\input xy
\usepackage[all]{xy}
\xyoption{line}
\xyoption{arrow}
\xyoption{color}
\SelectTips{cm}{}

\usepackage{tikz}
\usetikzlibrary{shapes,snakes}
\usetikzlibrary{decorations.markings}
\usetikzlibrary{decorations.pathreplacing}
\tikzstyle directed=[postaction={decorate,decoration={markings,
    mark=at position #1 with {\arrow{>}}}}]

\newcommand{\hackcenter}[1]{
 \xy (0,0)*{#1}; \endxy}

\usepackage{tikz-cd}
\tikzset{->-/.style={decoration={
  markings,
  mark=at position #1 with {\arrow{>}}},postaction={decorate}}}

\tikzset{middlearrow/.style={
        decoration={markings,
            mark= at position 0.5 with {\arrow{#1}} ,
        },
        postaction={decorate}
    }
}

\usepackage{bm}

\def\u0{{\underline{0}}}

\def\Q{\mathsf{Q}}

\theoremstyle{plain}
\newtheorem{theorem}{Theorem}

\newtheorem{lemma}[theorem]{Lemma}

\newtheorem{claim}[theorem]{Claim}

\newtheorem{notation}[theorem]{Notation}

\theoremstyle{definition}
\newtheorem{example}[theorem]{Example}
\newtheorem{definition}[theorem]{Definition}

\theoremstyle{definition}
\newtheorem{remark}[theorem]{Remark}

\numberwithin{equation}{section}
\numberwithin{theorem}{section}

\newcommand{\refequal}[1]{\xy {\ar@{=}^{#1}
(-1,0)*{};(1,0)*{}};
\endxy}

\newcommand{\cat}[1]{\ensuremath{\mbox{\bfseries {\upshape {#1}}}}}

\newcommand{\Hom}{{\rm Hom}}

\renewcommand{\to}{\rightarrow}

\def\Id{\mathrm{Id}}

\def\mf{\mathfrak}

\def\Q{\mathbb{Q}}

\numberwithin{equation}{section}
\def\ME#1{\textcolor[rgb]{0.40,0.00,0.90}{[ME: #1]}}%

\let\hat=\widehat

\let\epsilon=\varepsilon

\usepackage{bbm}
\def\C{{\mathbb{C}}}
\def\N{{\mathbbm N}}

\def\Z{{\mathbbm Z}}

\def\1{\mathbbm{1}}%
\usepackage[normalem]{ulem}

\newcommand\nc{\newcommand}
\nc\rnc{\renewcommand}
\nc\Kar{\operatorname{Kar}}
\nc\End{\operatorname{End}}

\nc\Sym{\operatorname{Sym}}

\allowdisplaybreaks

\title{A new presentation of the $\mf{osp}(1|2)$-polynomial link invariant and categorification}

\begin{document}
\setcounter{tocdepth}{3}

\author{Mark Ebert}
\email{markeber@usc.edu}
\address{Department of Mathematics\\ University of Southern California \\ Los Angeles, CA}
\maketitle
\begin{abstract}
There is a known connection between the $\mf{osp}(1|2n)$ polynomial knot invariant $J_K^n$ and the $\mf{so}(2n+1)$ knot invariant ${}_{so}J_K^n$ studied by Clark \cite{C} and Blumen \cite{Blu09}.
In the rank one case, the uncolored $U_{q}(\mf{osp}(1|2))$ link invariant is equal to the $U_{t^{-1}q}(\mf{sl}_2)$ link invariant where $t^2=-1$.
We define a skein relation similar to the Kauffman bracket, and use that to recover an oriented link invariant which coincides with Clark's uncolored $\mf{osp}(1|2)$-link invariant.
This definition also comes from the representation theory of $U_{q,\pi}(\mf{sl}_2)$, but using different methods from Clark.
We show that our invariant is easily categorified by a slightly modified version of Khovanov homology equipped with an extra $\Z_4$-grading.
We also construct a similarly modified version of Putyra's covering Khovanov homology \cite{P}.
This suggests that the similarity between the two invariants holds at the categorified level as well.
\end{abstract}

\section{Introduction}

The Reshetikhin-Turaev invariants are a family of link-invariants constructed from the representation theory of quantum groups \cite{RT90}.
Examples of such invariants include the HOMFLY-PT polynomial and the colored Jones polynomial.
Following Khovanov's pioneering work categorifying the Jones polynomial \cite{K00}, and the subsequent work of many others\cite{Kh03, Sus07, KR08a, KR08b, MS09}, Webster categorified these invariants \cite{Web} using the theory of categorified quantum groups developed by Khovanov, Lauda, and Rouquier \cite{Rou2,KL10}.
The study of these link homology theories has produced a myriad of results in many different fields.

A major development in the study of categorified link invariants was an alternative categorification of the Jones polynomial called odd Khovanov homology introduced by Ozsv\'ath, Rasmussen  and Szab\'o \cite{ORS}.
Odd Khovanov homology coincides with Khovanov homology when taken with $\Z_2$-coefficients, but is distinct from Khovanov homology in general \cite{Shum}.
In order to obtain a representation theoretic undestanding of this odd Khovanov homology, Ellis, Khovanov, and Lauda started a program studying odd analogues of structures related to $U_q(\mf{sl}_2)$ and its categorification \cite{EKL, EK,Lau-odd}.

These investigations into odd categorification turned out to be closely connected with independent parallel investigations into  Kac-Moody superalgebra categorifications~\cite{KKT, KKO, KKO2}, with the odd categorification of $\mf{sl}_2$ lifting the rank one Kac-Moody superalgebra.
These odd categorifications give categorifications of the theory of covering Kac-Moody algebras~\cite{HillWang,ClarkWang, CHW,CHW2}.
Covering algebras $U_{q,\pi}(\mf{g})$ generalize quantum enveloping algebras, depending on an additional parameter $\pi$ with $\pi^2 = 1$.  When $\pi=1$, it reduces to the usual quantum enveloping algebra $U_{q}(\mf{g})$, while the $\pi = -1$ specialization recovers the quantum group of a super Kac-Moody algebra.
The quantum algebra and quantum superalgebra can be identified by an automorphism of (an extension of) the covering quantum group called a twistor, which sends $\pi \to -\pi$ and $q\to t^{-1}q$ where $t^2=-1$ \cite{CFLW}.

In the rank one case, the $\pi=1$ specialization is $U_{q}(\mf{sl}_2)$, while for $\pi=-1$ it gives the quantum supergroup $U_{q}(\mf{osp}(1|2))$ associated with the lie superalgebra $\mf{osp}(1|2)$.
Ellis and Lauda defined a $2$-supercategory $\mf{U}$ that categorifies the rank one covering algebra $U_{q,\pi}(\mf{sl}_2)$ \cite{Lau-odd}.
This categorification was later greatly simplified and generalized by Brundan and Ellis~\cite{BE2}, where the 2-supercategory formalism was better developed, building off of their work~\cite{BE1}. 
However, a understanding of odd Khovanov homology in terms of higher representation theory of $\mf{U}$ still eludes us.

There is heurisitic evidence for the existence of odd link homology theories associated to $\mf{so}_{2n+1}$ from Mikhaylov and Witten \cite{Witten}.
If we are somehow able to develop an odd version of Webster's results in the future, then we can get odd link homology theories categorifying invariants obtained through the RT process on quantum supergroups.

The $U_{q}(\mf{osp}(1|2n))$ link invariant has been defined using Markov traces on the Birman-Wenzl-Murakami algebra\cite{Blu09}.
Clark \cite{C} produced the colored $U_q(\mf{osp}(1|2n))$ invariant from the representation theory of the quantum covering group $U_{q,\pi}(\mf{osp}(1|2n))$ a la Turaev \cite{Tu}.
A notable result from Clark's paper is that the $\mf{osp}(1|2n)$ invariant and the $\mf{so}(2n+1)$ invariant are basically the same using the theory of twistors.
If we represent a tangle as a composition $S$ of $U_{q,\pi}$-module homomorphisms, then Clark showed that twistors commute with $S$ up to an integral power of $\tau$.
\begin{theorem}{\cite[Theorem 4.24]{C}}
Let $K$ be any oriented knot, and let ${}_{\mf{osp}}J_K^\lambda(q), {}_{\mf{so}}J_K^\lambda(q)\in \Q(q,t)^\tau$ be the $\lambda$-colored $\mf{osp}$ and $\mf{so}$ knot invariants.
Then
  \[
  {}_{\mf{osp}}J_K^\lambda(q)=t^{\star(K,\lambda)}{}_{\mf{so}}J_K^\lambda(t^{-1}q)
  \]
  for some $\star(K,\lambda)\in \Z$.
\end{theorem}
In particular, the $U_q(\mf{osp}(1|2))$-invariant in the variable $q$ is just the Jones polynomial in the variable $\tau^{-1}q$.

If $T_+, T_-, T_0$ are unframed oriented links that differ as shown below
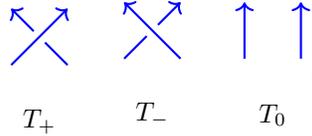
\begin{align*}
  \hackcenter{ \begin{tikzpicture} [scale=.75]
\draw[thick,   blue, ->] (0,0) to (1,1);
\draw[thick,   blue  ] (1,0) to (.6,.4);
\draw[thick,   blue, <-] (0,1) to (.4,.6);
 \node at (.5,-1) {$T_+$};
\end{tikzpicture}}
\qquad
  \hackcenter{ \begin{tikzpicture} [xscale=-1, scale=.75]
\draw[thick,   blue, ->] (0,0) to (1,1);
\draw[thick,   blue  ] (1,0) to (.6,.4);
\draw[thick,   blue, <-] (0,1) to (.4,.6);
 \node at (.5,-1) {$T_-$};
\end{tikzpicture}}
\qquad
  \hackcenter{ \begin{tikzpicture} [xscale=-1, scale=.75]
\draw[thick,   blue, ->] (0,0) to (0,1);
\draw[thick,   blue, ->] (1,0) to (1,1);
 \node at (.5,-1) {$T_0$};
\end{tikzpicture}},
\end{align*}
then
\begin{align}\label{osp-Skein}
  (\tau^{-1}q)^{-2} J_{T_+}^1-(\tau^{-1} q)^2 J_{T_-}^1=(\tau q^{-1}-\tau^{-1}q)J_{T_0}^1.
\end{align}
However, a categorification of this invariant has not been discovered as of yet.

\subsection{Results}
In this paper, we also construct an invariant $\hat{J}$ of oriented links using the representation theory of $\mf{osp}(1|2)$.
We first define a new Kauffman bracket-type Skein relation satisfying
\begin{align*}
   \left\langle
\hackcenter{ \begin{tikzpicture} [scale=.75]
\draw[thick,   blue] (0,0) to (1,1);
\draw[thick,   blue] (1,0) to (.6,.4);
\draw[thick,   blue] (0,1) to (.4,.6);
\end{tikzpicture}}
\right\rangle
\;\;&=\;\;
\tau \left\langle
\hackcenter{ \begin{tikzpicture} [scale=.75]
\draw[thick,   blue] (0,0) to (0,1);
\draw[thick,   blue] (1,0) to (1,1);
\end{tikzpicture}}
\right\rangle
-q\left\langle
\hackcenter{ \begin{tikzpicture} [scale=.75]
\draw[thick,   blue] (0,0).. controls ++(.25,.5) and ++(-.25,.5).. (1,0);
\draw[thick,   blue] (0,1).. controls ++(.25,-.5) and ++(-.25,-.5).. (1,1);
\end{tikzpicture}}
\right\rangle
\end{align*}
Then $\hat{J}$ is given by a renormalization of the bracket, and satisfies \ref{osp-Skein}.
Therefore $\hat{J}$ agrees with Clark's $\mf{osp}(1|2)$ invariant $J_L^1$ on any oriented link $L$.
We then show that the $\mf{osp}(1|2)$ link invariant formulated this way can be easily categorified with a slightly modified version of Khovanov homology with an additional $\Z_4$-grading.
We also construct a slight modification of Putyra's covering Khovanov homology \cite{P} that has an additional $\Z_4$ grading.
The idea of adding an additional grading to Khovanov homology has been done\cite{Man1}, but the connection to $\mf{osp}(1|2)$ is new.



\subsection{Acknowledgements}
I would like to thank Aaron Lauda, Joshua Sussan, and Pedro Vaz for their helpful comments on this article.
M.E was partially supported by NSF grants DMS-1902092 and DMS-2200419, the Army Research Office award W911NF-20-1-0075, and the Simons Foundation collaboration grant on New Structures in Low-dimensional topology.

\section{$\mf{osp}(1|2)$ braid invariant}
In this paper, $\tau$ is an indeterminant such that $\tau^4=1$, and $\pi:=\tau^2$. We also define the rings
\begin{align*}
  R^\tau &=R[\tau]/(\tau^4-1) &  R^\pi &=R[\pi]/(\pi^2-1)
\end{align*}

\begin{definition}
  Given $n\in \N$ and $d\in \Z^\tau[q,q^{-1}]$, the Temperley-Lieb algebra $TL_n(d)$ is defined as the $\Z^\tau[q,q^{-1}]$-algebra with generators $\{1,E_1, \dots, E_{n-1}\}$ with relations.
  \begin{enumerate}
    \item $E_iE_j=E_jE_i$ for $|i-j|>1$
    \item $E_{i}E_{i\pm 1}E_i=E_i$
    \item $E_i^2=dE_i$
  \end{enumerate}
  We can also think of it as a purely even superalgebra.
\end{definition}

For this section, set
\[
d:=\tau^{-1}q+ \tau q^{-1}\in \Z^{\tau}[q,q^{-1}]
\]

\begin{definition}
  The braid group $B_n$ on $n$ strands is the group
  \[B_n= \langle \sigma_1, \dots, \sigma_{n-1} \mid \sigma_{i}\sigma_{i+1}\sigma_{i}=\sigma_{i+1}\sigma_{i}\sigma_{i+1}, \; \sigma_i\sigma_j = \sigma_j \sigma_i \rangle\]
  where $|i-j|>1$.
\end{definition}

We can represent elements of the braid group and the Tempereley Lieb algebra diagrammatically
\begin{align*}
\sigma_i=
\hackcenter{ \begin{tikzpicture} [scale=.75]
\draw[thick,   blue, ->] (0,0) to (1,1);
\draw[thick,   blue  ] (1,0) to (.6,.4);
\draw[thick,   blue, <-] (0,1) to (.4,.6);
\end{tikzpicture}}& \qquad
\sigma_i^{-1}=
\hackcenter{ \begin{tikzpicture} [xscale=-1,scale=.75]
\draw[thick,   blue, ->] (0,0) to (1,1);
\draw[thick,   blue] (1,0) to (.6,.4);
\draw[thick,   blue, <-] (0,1) to (.4,.6);
\end{tikzpicture}} & \qquad
  Id:=
\hackcenter{ \begin{tikzpicture} [scale=.75]
\draw[thick,   blue] (0,0) to (0,1);
\draw[thick,   blue] (1,0) to (1,1);
\end{tikzpicture}} & \qquad E_i:=
\hackcenter{ \begin{tikzpicture} [scale=.75]
\draw[thick,   blue] (0,0).. controls ++(.25,.5) and ++(-.25,.5).. (1,0);
\draw[thick,   blue] (0,1).. controls ++(.25,-.5) and ++(-.25,-.5).. (1,1);
\end{tikzpicture}}
\end{align*}
We take braids to be oriented upwards, and so remove the arrows from now on.

\begin{definition}
  Given $B\in B_n$ thought of as an upward oriented braid, we define the trace closure $B^{tr}$ as the link formed by closing off the braid to the right. 
\end{definition}
For example, 
\begin{align*}
  \sigma =  \hackcenter{ \begin{tikzpicture} [scale=.75]
\draw[thick,   blue, -] (0,0) to (1,1);
\draw[thick,   blue  ] (1,0) to (.6,.4);
\draw[thick,   blue, -] (0,1) to (.4,.6);
\end{tikzpicture}} \;\; \in B_2& \qquad \text{implies}\qquad
\sigma^{tr}=
  \hackcenter{\begin{tikzpicture} [scale=.75]
\draw[thick,   blue] (0,0) to (1,1);
\draw[thick,   blue] (1,0) to (.6,.4);
\draw[thick,   blue] (0,1) to (.4,.6);
\draw[thick, blue] (1,1) .. controls ++(.5,.25) and ++(.5, -.25) ..(1,0);
\draw[thick, blue] (0,1) .. controls ++(-.25,.25) and ++(-.2, .25) ..(1.5,1.5);
\draw[thick, blue] (0,0) .. controls ++(-.25,-.25) and ++(-.2, -.25) ..(1.5,-0.5);
\draw[thick, blue] (1.5,1.5) .. controls++ (0.5,-0.5) and ++(0.5,0.5) ..(1.5,-0.5);
\end{tikzpicture}}.
\end{align*}

\begin{definition}
  Given a braid $B\in B_n$, define the bracket polynomial $\langle B \rangle\in \Z^\tau[q,q^{-1}]$ by $\langle B \rangle := \langle B^{tr} \rangle$, where $\langle B^{tr} \rangle$ is characterized by the following local relations on the braid closure.
\end{definition}

\begin{align}
  \left\langle
\hackcenter{ \begin{tikzpicture} [scale=.75]
\draw[thick,   blue] (0,0) to (1,1);
\draw[thick,   blue] (1,0) to (.6,.4);
\draw[thick,   blue] (0,1) to (.4,.6);
\end{tikzpicture}}
\right\rangle
\;\;&=\;\;
\tau \left\langle
\hackcenter{ \begin{tikzpicture} [scale=.75]
\draw[thick,   blue] (0,0) to (0,1);
\draw[thick,   blue] (1,0) to (1,1);
\end{tikzpicture}}
\right\rangle
-q\left\langle
\hackcenter{ \begin{tikzpicture} [scale=.75]
\draw[thick,   blue] (0,0).. controls ++(.25,.5) and ++(-.25,.5).. (1,0);
\draw[thick,   blue] (0,1).. controls ++(.25,-.5) and ++(-.25,-.5).. (1,1);
\end{tikzpicture}}
\right\rangle
\qquad\text{OR}\;\;
    \left\langle
\hackcenter{ \begin{tikzpicture} [xscale=-1,scale=.75]
\draw[thick,   blue] (0,0) to (1,1);
\draw[thick,   blue] (1,0) to (.6,.4);
\draw[thick,   blue] (0,1) to (.4,.6);
\end{tikzpicture}}
\right\rangle
\;\;=\;\;
\tau^3 \left\langle
\hackcenter{ \begin{tikzpicture} [scale=.75]
\draw[thick,   blue] (0,0) to (0,1);
\draw[thick,   blue] (1,0) to (1,1);
\end{tikzpicture}}
\right\rangle
-q^{-1} \left\langle
\hackcenter{ \begin{tikzpicture} [scale=.75]
\draw[thick,   blue] (0,0).. controls ++(.25,.5) and ++(-.25,.5).. (1,0);
\draw[thick,   blue] (0,1).. controls ++(.25,-.5) and ++(-.25,-.5).. (1,1);
\end{tikzpicture}}\right\rangle
\label{bracket-poly2-1}
\\
\langle B_1 \sqcup B_2 \rangle &:= \langle B_1 \rangle \langle B_2 \rangle
\qquad
 \left\langle
 \hackcenter{ \begin{tikzpicture} [scale=.75]
\draw[thick,   blue] (0,0).. controls ++(.15,.5) and ++(-.15,.5).. (1,0);
\draw[thick,   blue] (0,0).. controls ++(.15,-.5) and ++(-.15,-.5).. (1,0);
\end{tikzpicture}}
\right\rangle
\;\;=\;\;
d
 \qquad \langle \emptyset \rangle=1
 \label{bracket-poly2-2}
\\
\pi \;\;
\hackcenter{ \begin{tikzpicture} [ scale=.5]
%
\draw[thick,   blue] (0,0).. controls ++(.25,.5) and ++(-.25,.5).. (1,0);
\draw[thick,   blue] (1,0).. controls ++(.25,-.5) and ++(-.25,-.5).. (2,0);
\draw[thick, blue] (2,0) to (2,1);
\draw[thick, blue] (0,0) to (0,-1);
\end{tikzpicture}}
\;\;&=\;\;
-\tau \;\;
\hackcenter{ \begin{tikzpicture} [scale=.5]
\draw[thick,   blue] (0,-1) to (0,2);
\end{tikzpicture}}
\;\;=\;\;
  \hackcenter{ \begin{tikzpicture} [xscale=-1, scale=.5]
%
\draw[thick,   blue] (0,0).. controls ++(.25,.5) and ++(-.25,.5).. (1,0);
\draw[thick,   blue] (1,0).. controls ++(.25,-.5) and ++(-.25,-.5).. (2,0);
\draw[thick, blue] (2,0) to (2,1);
\draw[thick, blue] (0,0) to (0,-1);
\end{tikzpicture}}
\qquad
\text{WHICH IMPLIES}
\qquad
  \hackcenter{ \begin{tikzpicture} [scale=.5]
\draw[thick,   blue] (0,1).. controls ++(.25,-.5) and ++(-.25,-.5).. (1,1);
\draw[thick,   blue] (1,1).. controls ++(.25,.5) and ++(-.25,.5).. (2,1);
\draw[thick,   blue] (0,0).. controls ++(.25,.5) and ++(-.25,.5).. (1,0);
\draw[thick,   blue] (1,0).. controls ++(.25,-.5) and ++(-.25,-.5).. (2,0);
\draw[thick, blue] (2,0) to (2,1);
\draw[thick, blue] (0,0) to (0,-1);
\draw[thick, blue] (0,1) to (0,2);
\end{tikzpicture}}
\;\;=\;\;
  \hackcenter{ \begin{tikzpicture} [xscale=-1, scale=.5]
\draw[thick,   blue] (0,1).. controls ++(.25,-.5) and ++(-.25,-.5).. (1,1);
\draw[thick,   blue] (1,1).. controls ++(.25,.5) and ++(-.25,.5).. (2,1);
\draw[thick,   blue] (0,0).. controls ++(.25,.5) and ++(-.25,.5).. (1,0);
\draw[thick,   blue] (1,0).. controls ++(.25,-.5) and ++(-.25,-.5).. (2,0);
\draw[thick, blue] (2,0) to (2,1);
\draw[thick, blue] (0,0) to (0,-1);
\draw[thick, blue] (0,1) to (0,2);
\end{tikzpicture}}
\;\;=\;\;
\hackcenter{ \begin{tikzpicture} [scale=.5]
\draw[thick,   blue] (0,-1) to (0,2);
\end{tikzpicture}}
\label{bracket-poly2-3}
\end{align}

\begin{remark}
  The second relation of \eqref{bracket-poly2-1} follows from the first and vice versa since the braid group relations imply that $\langle \sigma_i \sigma_i^{-1}\rangle=\langle Id \rangle $. We include both here for convenience.
\end{remark}

\begin{remark}\label{rmk:closed-loops}
  The first relation of \eqref{bracket-poly2-3} clearly implies second relation, and it does hold in the representation theoretic perspective.
  However, we don't actually need it to define the bracket polynomial.
  We can in fact reduce all closed loops to circles using the second relation. The second relation provides an easy connection to the Temperley Lieb algebra.
\end{remark}

\begin{definition}\label{def:state}
A state $z$ of a braid $B$ to be a choice of resolution
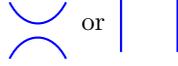
\begin{align*}
  \hackcenter{ \begin{tikzpicture} [scale=.75]
\draw[thick,   blue] (0,0).. controls ++(.25,.5) and ++(-.25,.5).. (1,0);
\draw[thick,   blue] (0,1).. controls ++(.25,-.5) and ++(-.25,-.5).. (1,1);
\end{tikzpicture}}
\;\; \text{or} \;\;
\hackcenter{ \begin{tikzpicture} [scale=.75]
\draw[thick,   blue] (0,0) to (0,1);
\draw[thick,   blue] (1,0) to (1,1);
\end{tikzpicture}}
\end{align*}
at each crossing of $B$.
\end{definition}

We call each resolution of a crossing a $0$-resolution or a $1$-resolution according to the figure below.
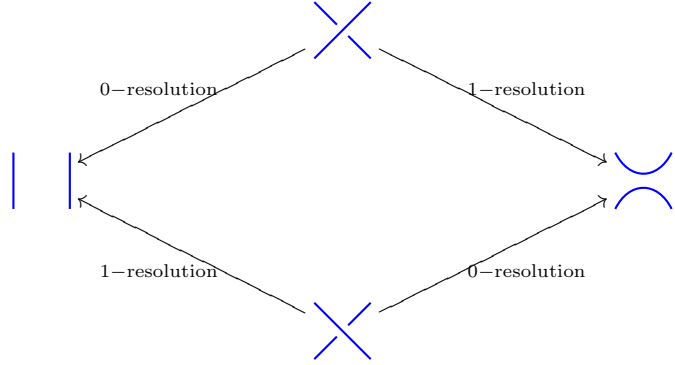
\[
\xy
  (0,20)*+{ \hackcenter{ \begin{tikzpicture} [scale=.75]
\draw[thick,   blue] (0,0) to (1,1);
\draw[thick,   blue] (1,0) to (.6,.4);
\draw[thick,   blue] (0,1) to (.4,.6);
\end{tikzpicture}}}="T";
  (0,-20)*+{ \hackcenter{ \begin{tikzpicture} [xscale=-1,scale=.75]
\draw[thick,   blue] (0,0) to (1,1);
\draw[thick,   blue] (1,0) to (.6,.4);
\draw[thick,   blue] (0,1) to (.4,.6);
\end{tikzpicture}}}="B";
  (-40,0)*+{ \hackcenter{ \begin{tikzpicture} [scale=.75]
\draw[thick,   blue] (0,0) to (0,1);
\draw[thick,   blue] (1,0) to (1,1);
\end{tikzpicture}}
 }="L";
  (40,0)*+{\hackcenter{ \begin{tikzpicture} [scale=.75]
\draw[thick,   blue] (0,0).. controls ++(.25,.5) and ++(-.25,.5).. (1,0);
\draw[thick,   blue] (0,1).. controls ++(.25,-.5) and ++(-.25,-.5).. (1,1);
\end{tikzpicture}} }="R";
  {\ar_-{0-\text{resolution}} "T";"L"};
  {\ar^-{1-\text{resolution}} "T";"R"};
  {\ar^-{1-\text{resolution}} "B";"L"};
  {\ar_-{0-\text{resolution}} "B";"R"};
\endxy
\]

If $B$ has $m$ crossings (which we label), then states of $B$ correspond to sequences $z\in \{0,1\}^m$ where the $i$th element says how we resolved the ith crossing of $B$. We can define the cube of resolutions as follows.
Take the graph with vertices given by sequences $z\in \{0,1\}^m$ and edges $z\to w$ where $w$ can be obtained from $z$ by changing exactly one $0$ in $z$ to a $1$.
To each vertex $z$, associate a diagram $B_z'$ given by resolving each crossing according to $z$ (remember not to evaluate the bubbles).
This is known as the cube of resolutions, and also is defined for knots and links.
When we discuss categorification later in section~\ref{section:cat}, we will associate a saddle cobordism $(B_{z}')^{tr}\to (B_{w}')^{tr}$ to each edge $z\to w$.

\begin{definition}
  For a braid $B\in B_n$ and a state $z$, we define a polynomial $f_z(q,\tau)\in \Z^\tau[q,q^{-1}]$ and a diagram $B_z\in TL_n(d)$ by first replacing each crossing of $B$ with the term in \eqref{bracket-poly2-1} containing the resolution specified by $z$ and then evaluating each closed loop using \eqref{bracket-poly2-2} and \eqref{bracket-poly2-3}.
  The result is an element in $TL_n(d)$ given by a diagram in $TL_n(d)$ multiplied by a polynomial in $\Z^\tau[q, q^{-1}]$.
  We call the resulting diagram $B_z$, and the resulting polynomial $f_{z}(q,\tau)$.
\end{definition}

It is not difficult to see that $B_z$ is the diagram obtained by deleting all of the closed loops in $B_z'$.
We can write $f_z(q,\tau)$ in the form
\[f_z(q,\tau)=\pm \tau^{z_0}q^{z_1} d^k,\]
where $z_0,z_1$ are determined by \eqref{bracket-poly2-1}, and $k$ is the number of closed loops that were removed from $B_z'$ to get $B_z$.

It is evident from our construction that the following is true.
\begin{align}\label{state-sum}
  \langle B \rangle=\sum_z f_z(q,\tau) \langle B_z^{tr}\rangle
\end{align}
where $\langle B_z^{tr} \rangle:=\langle (B_z)^{tr} \rangle$ is evaluated using \eqref{bracket-poly2-2} and \eqref{bracket-poly2-3}.


In the same way that one renormalizes the Kauffman bracket to get the Jones polynomial, we also renormalize $\langle - \rangle$ to get an invariant.

\begin{definition}\label{writhe}
  The writhe of a braid (or link) $T$ is defined by setting
\begin{align*}
wr\left(
\hackcenter{ \begin{tikzpicture} [scale=.75]
\draw[thick,   blue] (0,0) to (1,1);
\draw[thick,   blue] (1,0) to (.6,.4);
\draw[thick,   blue] (0,1) to (.4,.6);
\end{tikzpicture}}
\right)=1 & \qquad
wr\left(
\hackcenter{ \begin{tikzpicture} [xscale=-1,scale=.75]
\draw[thick,   blue] (0,0) to (1,1);
\draw[thick,   blue] (1,0) to (.6,.4);
\draw[thick,   blue] (0,1) to (.4,.6);
\end{tikzpicture}}
\right)=-1
& \qquad
wr(T)=\sum wr(X)
\end{align*}
where the sum is over all crossings $X$ in $T$.
We often write this as $wr(T)=n_+-n_-$, where $n_+$ and $n_-$ denote the number of positive and negative crossings in $T$ respectively.
\end{definition}

We are now ready to define the $\mf{osp}(1|2)$-polynomial invariant on a braid $T$. This will turn out define a invariant of oriented links.
\begin{definition}
  The uncolored $\mf{osp}(1|2)$-polynomial $\hat{J}$ is defined as a renormalization of the bracket polynomial
  \footnote{This new definition $\hat{J}(-)$ of the $\mf{osp}(1|2)$-polynomial is designed to agree with the uncolored $\mf{osp}(1|2)$-polynomial $J_{-}^l$ given in \cite[Example 3.10]{C}. The main difference between the two definitions is the type of skein relation used as mentioned in the introduction. The skein relation he uses can be derived from the skein relation I use, but not the other way around.}
  .
  \[
  \hat{J}(T)=(\tau^2 q)^{wr(T)}\langle T \rangle
  \]
\end{definition}

To extend our discussion to oriented links, we recall Alexander's theorem, which states any oriented link can be obtained as the trace closure of a braid.
It is possible that two different braids can have the same trace closure, so $\hat{J}$ extends to an invariant of oriented links if and only if it is well defined on isomorphism classes of closures of braids.
If this is the case, then if $L$ can be realized as $T^{tr}$ for some braid $T$, then $\hat{J}(L)=\hat{J}(T^{tr})=\hat{J}(T)$.

We say that $\hat{J}$ preserves the Markov moves if it satisfies the conditions M1 and M2 listed below.
\begin{enumerate}
  \item[M1] $\hat{J}(T\sigma_n^{\pm 1})=\hat{J}(T)$ for all $T\in B_n$
  \item[M2] $\hat{J}(AB)=\hat{J}(BA)$ for all $A,B \in B_n$
\end{enumerate}
If $\hat{J}$ preserves the Markov moves as above, then Markov's theorem tells us that $\hat{J}$ is well defined on isomorphism classes of closures of braids, and is thus an invariant of oriented links.
We will now show that $\hat{J}$ preserves the Markov moves.

\begin{lemma}\label{M1-pf}
  $\hat{J}$ satisfies M1.
\end{lemma}

\begin{proof}
Assume for the duration of the proof that $T\in B_n$.
We need to show that $\hat{J}(T)=\hat{J}(T\sigma_{n}^{\pm 1})$.
We only show the proof for $T\sigma_n^{+1}$, as the proof of the other case is essentially the same.

$T\sigma_n \in B_{n+1}$ has the following diagram
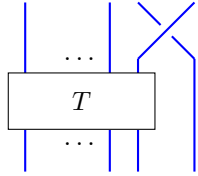
\[
\hackcenter{\begin{tikzpicture} [scale=0.75]
\draw[thick,  blue] (.5,-.5) to (.5,2.5);
\draw[thick,  blue] (-1,-.5) to (-1,2.5);
\draw[thick,  blue] (1,-.5) to (1,1.5);
\draw[fill=white!20,] (-1.3,.25) rectangle (1.3,1.25);
\node at (0,0) {$\dots$};
\node at (0,1.5) {$\dots$};
 \node at (0,.75) {$T$};
 \draw[thick,   blue] (1,1.5) to (2,2.5);
\draw[thick,   blue] (2,1.5) to (1.6,1.9);
\draw[thick,   blue] (1,2.5) to (1.4,2.1);
\draw[thick,blue] (2,1.5) to (2,-.5);
\end{tikzpicture}}
\]
%
%
Using \eqref{bracket-poly2-1}, we see that the trace closure of this diagram is equal to the trace closure of
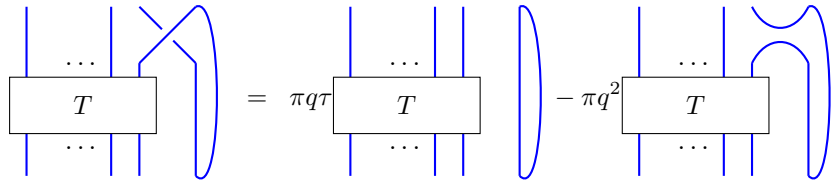
\[
\hackcenter{\begin{tikzpicture} [scale=0.75]
\draw[thick,  blue] (.5,-.5) to (.5,2.5);
\draw[thick,  blue] (-1,-.5) to (-1,2.5);
\draw[thick,  blue] (1,-.5) to (1,1.5);
\draw[fill=white!20,] (-1.3,.25) rectangle (1.3,1.25);
\node at (0,0) {$\dots$};
\node at (0,1.5) {$\dots$};
 \node at (0,.75) {$T$};
 \draw[thick,   blue] (1,1.5) to (2,2.5);
\draw[thick,   blue] (2,1.5) to (1.6,1.9);
\draw[thick,   blue] (1,2.5) to (1.4,2.1);
\draw[thick,blue] (2,1.5) to (2,-.5);
\draw[thick, blue] (2,2.5) .. controls ++(.5,.25) and ++(.5, -.5) ..(2,-.5);
\end{tikzpicture}}
\;\;=\;\;\pi q
\tau
\hackcenter{\begin{tikzpicture} [scale=0.75]
\draw[thick,  blue] (.5,-.5) to (.5,2.5);
\draw[thick,  blue] (-1,-.5) to (-1,2.5);
\draw[thick,  blue] (1,-.5) to (1,2.5);
\draw[fill=white!20,] (-1.3,.25) rectangle (1.3,1.25);
\node at (0,0) {$\dots$};
\node at (0,1.5) {$\dots$};
 \node at (0,.75) {$T$};
\draw[thick,blue] (2,2.5) to (2,-.5);
\draw[thick, blue] (2,2.5) .. controls ++(.5,.25) and ++(.5, -.5) ..(2,-.5);
\end{tikzpicture}}
-\pi q^2
\hackcenter{\begin{tikzpicture} [scale=0.75]
\draw[thick,  blue] (.5,-.5) to (.5,2.5);
\draw[thick,  blue] (-1,-.5) to (-1,2.5);
\draw[thick,  blue] (1,-.5) to (1,1.5);
\draw[fill=white!20,] (-1.3,.25) rectangle (1.3,1.25);
\node at (0,0) {$\dots$};
\node at (0,1.5) {$\dots$};
 \node at (0,.75) {$T$};
\draw[thick,   blue] (1,1.5).. controls ++(.25,.5) and ++(-.25,.5).. (2,1.5);
\draw[thick,   blue] (1,2.5).. controls ++(.25,-.5) and ++(-.25,-.5).. (2,2.5);
\draw[thick,blue] (2,1.5) to (2,-.5);
\draw[thick, blue] (2,2.5) .. controls ++(.5,.25) and ++(.5, -.5) ..(2,-.5);
\end{tikzpicture}}
\]

Hence, $\hat{J}(T\sigma_n)=(\tau^3q(\tau^3 q+\tau q^{-1})-\tau^2q^2)\hat{J}(T)=\hat{J}(T)$.
\end{proof}
\begin{lemma}\label{M2-pf}
  $\hat{J}$ satisfies M2.
\end{lemma}
\begin{proof}
For $A,B \in B_n$ we need to show that $\hat{J}(AB)=\hat{J}(BA)$ for any $A,B\in B_n$.

Let $x$ be a state of $A$ and $y$ be a state of $B$. Then we get canonical states $xy$ and $yx$ of $AB$ and $BA$ respectively.
Let $\ell_{AB}$ be the number of loops in $A_xB_y$. Then $f_{xy}=f_xf_y d^{\ell_{AB}}$ and $\langle (A_xB_y)^{tr}\rangle =d^{\ell_{AB}}\langle ((AB)_{xy})^{tr} \rangle$ by construction.
Therefore, $f_{xy}\langle ((AB)_{xy})^{tr}\rangle=f_xf_y \langle (A_xB_y)^{tr}\rangle $.
Similarly, $f_{yx}\langle ((BA)_{yx})^{tr} \rangle=f_yf_x \langle (B_yA_x)^{tr}\rangle$.

Therefore, it suffices to prove that $\langle (A_xB_y)^{tr}\rangle=\langle (B_yA_x)^{tr} \rangle$ for an arbitrary choice of states $x$ of $A$ and $y$ of $B$.
Since $A_x, B_y\in TL_n(d)$, we can use properties of the Temperley Lieb algebra to prove this without much effort.
Let $tr':TL_n(d)\to \Z^{\tau}[q,q^{- 1}]$ denote the Markov trace on $TL_n(d)$ (see \cite[Definition 2.12]{AJL}).
Then, $\langle (A_xB_y)^{tr} \rangle=d^n tr'(A_xB_y)$ by definition.
The Markov trace satisfies $tr'(AB)=tr'(BA)$, hence we have that $\langle (A_xB_y)^{tr} \rangle= \langle (B_yA_x)^{tr} \rangle$ for any states $x,y$.
Therefore,
\[
\hat{J}(AB)=\sum_{xy} f_{xy}\langle ((AB)_{xy})^{tr}\rangle= \sum_{x,y}f_xf_y\langle (A_xB_y)^{tr}\rangle=\sum_{x,y}f_xf_y \langle (B_yA_x)^{tr}\rangle=\sum_{yx} f_{yx}\langle ((BA)_{yx})^{tr}\rangle=\hat{J}(BA)
\]
\end{proof}

Hence, we have proven the following theorem.
\begin{theorem}\label{thm:decat-braid-invar}
  $\hat{J}$ is an invariant of oriented links.
\end{theorem}

\section{Explanation using representation theory}

Blumen\cite{Blu09} was the first to show the connection between the $\mf{so}$ and $\mf{osp}$ link invariants, but Clark\cite{C} was able to formulate the $\mf{osp}(1|2)$-invariant using the representation theory of the quantum covering group $U_{q,\pi}(\mf{sl}_2)$.
We give an overview of the quantum covering group $U_{q,\pi}$ defined by Clark Hill and Wang \cite{CHW} and its representation theory in section~\ref{sec:covering-group} and section~\ref{sec:U-modules} respectively.
We then review Clark's representation-theoretic definition of the invariant from \cite{C} in section~\ref{Sec:Clark-method}.
Finally, we give a representation-theoretic definition of our invariant $\hat{J}$ in section~\ref{Sec:Our-method}.

\subsection{Parameters}\label{sec:3.1}

\begin{notation}
  We borrow much of our notation in this section from \cite{C}. Given an element $m$ in some $\Z\times \Z_2$-graded module $M$, we denote $p(m)$ to be the parity of $m$, $|m|$ denotes the $\Z$-degree of $m$ and $||m||=(|m|,p(m))$.
\end{notation}

Let us review some notions from \cite[Section 2.2]{C}.

Let $t\in \C$ such that $t^2=-1$. Let $q$ be a formal parameter, let $\tau$ be a formal parameter such that $\tau^4=1$, and define $\pi:=\tau^2$.
Recall the notation $R^\pi=R[\pi]/(\pi^2-1)$ and $R^\tau = R[\tau]/(\tau^4-1)$ for a commutative ring $R$ with 1.
Our base ring is $\Q(q,t)^\tau$.

For $x\in \{\pm 1, \pm t\}$, let $\Q(q,t)_x:=\Q(q,t)$ viewed as a $\Q(q,t)^\tau$-module where $\tau$ acts by multiplication by $x$.
\begin{definition}
  The specialization of $\Q(q,t)^\tau$-module $M$ at $x$ is the $\Q(q,t)$-module \[M|_{\tau=x}=\Q(q,t)_x\otimes_{\Q(q,t)^\tau}M\]
\end{definition}

$\Q(q,t)^\tau$ has orthogonal idempotents
\begin{align*}
  \varepsilon_{t^k}&=\frac{1+(t^k\tau)+(t^k\tau)^2+(t^k\tau)^3}{4} & 0\leq k\leq 3,
\end{align*}
such that $\Q(q,t)^\tau=\bigoplus_{k\in \Z_4}Q(q,t)\varepsilon_{t^k}$.

Furthermore, since $\tau \varepsilon_x=x\varepsilon_x$, we see that
\[
M|_{\tau=x}=\varepsilon_x M.
\]

\subsection{Quantum covering group}\label{sec:covering-group}

\begin{definition}{\cite{CHW}}
  The quantum covering group $U_{q,\pi}:=U_{q,\pi}(\mf{sl}_2)$ is the $\Q(q,t)^\tau$-algebra with generators $E,F,K^{\pm 1},J$, subject to the relations
  \begin{align}\label{quantum-covering-group}
    J^2&=1 & JK&=KJ \\
    JE&=EJ & JF&=FJ \\
    KE&=q^2EK & KF&=q^{-2}FK \\
    EF-\pi FE &= \frac{JK-K^{-1}}{\pi q-q^{-1}}
  \end{align}
\end{definition}

It has structure of a Hopf superalgebra\cite{CHW,CHW2, C} but with braiding replaced by $x\otimes y \to \pi^{p(x)p(y)}y\otimes x$.
The coproduct $\Delta:U_{q,\pi} \to U_{q,\pi}\otimes U_{q,\pi} $ is
\begin{align*}
  \Delta(E)&=E\otimes 1+ JK\otimes E & \Delta(F)&=F\otimes K^{-1}+1\otimes F & \Delta(K)&=K\otimes K & \Delta(J)&=J\otimes J
\end{align*}

The antipode $S:U_{q,\pi}\to U_{q,\pi}$ is
\begin{align*}
  S(E)&=-J^{-1}K^{-1}E & S(F)&=-FK  & S(K)&=K^{-1} & S(J)&=J
\end{align*}

The counit $\varepsilon:U \to \Q(q,t)^{\tau}$
\begin{align*}
  \varepsilon(E)=\varepsilon(F)=0 \qquad \varepsilon(J)=\varepsilon(K)=1
\end{align*}

\subsection{$U_{q,\pi}$-modules}\label{sec:U-modules}

\begin{definition}{\cite[Section 2.4]{C}}
  A weight $U_{q,\pi}$-module $M$ is a $\Z\times \Z_2 $ graded $U$-module, whose grading is compatible with the grading on $U$ such that
  \begin{align*}
    M=\bigoplus_{\lambda \in \Z} M_{\lambda,0}\oplus M_{\lambda,1} &\qquad M_{\lambda,s}=\{ m\in M \mid p(m)=s, \; Jm=\pi^\lambda m, \;Km=q^\lambda m \}
  \end{align*}
\end{definition}

The simple weight $U$-module of highest weight $m$ with highest weight vector $v_m$ of even parity ($||v_m||=(m,0)$) is denoted $V(m)$.
It has basis $\{v_{m-2k}:=F^{(k)}v_m\}_{0\leq k\leq m}$. By convention we have $p(v_{m-2k})\equiv k (mod\; 2)$. The action of $U_{q,\pi}$ on $V(m)$ is given by
\begin{align*}
  Ev_{m-2k}&=\pi^{k-1}[m+1-k]v_{m-2(k-1)}  &\qquad Kv_{m-2k}&=q^{m-2k}v_{m-2k}\\
  Fv_{m-2k}&=[k+1]v_{m-2(k+1)} & \qquad Jv_{m-2k}&=\pi^{m-2k}v_{m-2k}
\end{align*}

Define the parity reversed module $\Pi M$ to be $M$ as a vector space equipped with the same action of $U_{q,\pi}$, but with $(\Pi M)_{\lambda, s}=M_{\lambda,1-s}$.
In particular, the identity map induces an odd $U$-module isomorphism $\gamma_{V(m)}:V(m) \to \Pi V(m)$.

For the purposes of this paper, it is enough to consider the case where $m=1$.
As a superspace, the $U_{q,\pi}$-module $V:=V(1)$ is given by the span of elements $\{v_+, v_-\}$ with $\Z\times \Z_2$ degrees $||v_+||=(1,0)$ and $||v_-||=(-1,1)$.
The action of $U_{q,\pi}$ on $V$ is given by
\begin{align*}
  Ev_{-}&=v_+  & Ev_+&=0 &\qquad Kv_{\pm}&=q^{\pm 1}v_{\pm}\\
  Fv_{+}&=v_{-} & Fv_-&=0 & \qquad Jv_{\pm}&=\pi^{\pm 1}v_{\pm}
\end{align*}

Then $\Pi V=\text{span} \{v_+', v_-'\}$ as a superspace where the $\Z\times \Z_2$ degrees of the basis elements are $||v_+'||=(1,1)$ and $||v_-'||=(-1,0)$.
By definition, the action of $U_{q,\pi}$ on $\Pi V$ is the same as the action on $V$.
\begin{align*}
  Ev_{-}'&=v_+'  & Ev_+'&=0 &\qquad Kv_{\pm}'&=q^{\pm 1}v_{\pm}'\\
  Fv_{+}'&=v_{-}' & Fv_-'&=0 & \qquad Jv_{\pm}'&=\pi^{\pm 1}v_{\pm}'
\end{align*}

Thus, the identity map induces an odd $U$-module isomorphism $\gamma_V: V \to \Pi V$ which sends $v_\pm \to v_\pm'$.
The dual module $M^*$ is defined to be
\begin{align*}
  M^*=\bigoplus_{\lambda \in \Z} (M_{\lambda,0})^*\oplus (M_{\lambda,1})^* &\qquad (M_{\lambda,s})^*=\Hom_{Q(q,t)^\tau}(M_{\lambda,s},Q(q,t)^\tau)
\end{align*}
The action on $M^*$ is induced by the Hopf superalgebra structure of $U_{q,\pi}$: For $x\in U_{q,\pi}$,
\begin{align*}
  xf(v)=\pi^{p(f)p(x)}f(S(x)v).
\end{align*}

Let us consider $M=V(m)$.
Then $V(m)^*$ has a dual basis $\{v_{m-2k}^*\}_{0 \leq k\leq m}$ with $||v_{m-2k}^*||=(2k-m,k)$.
\begin{align*}
  Ev_{m-2k}^* &= -(\pi q^{-1})^{m-2k}[m-k]v_{m-2(k+1)}^* \\
  Fv_{m-2k}^* &= -\pi^k q^{m-2k+2}[k]v_{m-2(k-1)}^*
\end{align*}

We prove the action of $F$ is correct.
\begin{align*}
  F v_{m-2k}^*(v_\alpha) &= \pi^k v_{m-2k}^*(S(F)v_{\alpha}) = \pi^k v_{m-2k}^*(-FKv_{\alpha}) \\
  &=-\pi^k\delta_{\alpha-2, m-2k}q^\alpha v_{m-2k}^*(Fv_{\alpha})= -\delta_{\alpha-2, m-2k} \pi^k q^{m-2k+2} v_{m-2k}^*(Fv_{m-2(k-1)}) \\
  &=-\delta_{\alpha-2, m-2k} \pi^k q^{m-2k+2} [k]v_{m-2k}^*(v_{\alpha-2})=-\pi^kq^{m-2k+2}[k] v_{m-2(k-1)}^*(v_\alpha)
\end{align*}

$V(m)^*$ is a weight module with highest weight vector $v_{-m}^*$. This has degree $||v_{-m}^*||=(m,m)$, so there is an even $U_{q,\pi}$-module isomorphism $\Pi^m V(m) \to V(m)^*$.

\subsection{Clark's $\mf{osp}(1|2)$ invariant}\label{Sec:Clark-method}
In this subsection, we review how Clark obtains the uncolored $\mf{osp}(1|2)$ link invariant in simplest case done in \cite[Example 3.10]{C}.
Clark actually defines a colored $\mf{osp}(1|2n)$-invariant, but we only consider the uncolored case with $n=1$.
All of the material in this subsection comes from \cite[Section 2, Section 3]{C}.

\begin{notation}
  We use the ordered basis
  \begin{align*}
    \{v_{+}\otimes v_{+}, v_{+}\otimes v_{-}, v_{-}\otimes v_{+}, v_{-}\otimes v_{-}\} &\qquad \text{for} \; V\otimes V \\
    \{v_{+}^*\otimes v_{+}, v_{+}^*\otimes v_{-}, v_{-}^*\otimes v_{+}, v_{-}^*\otimes v_{-}\} &\qquad \text{for} \; V^*\otimes V  \\
    \{v_{+}\otimes v_{+}^*, v_{+}\otimes v_{-}^*, v_{-}\otimes v_{+}^*, v_{-}\otimes v_{-}^*\} &\qquad \text{for} \; V\otimes V^*
  \end{align*}
\end{notation}

Recall the $\Q(q,t)^{\tau}$-linear maps from \cite[Lemma 2.7]{C}.
\begin{itemize}
  \item $\text{ev}:V(1)^* \otimes V(1) \to \Q(q,t)^{\tau}$ is the map $v^*\otimes w \to v^*(w)$. It is given by the matrix
  $\begin{pmatrix}
    1 & 0 & 0 & 1
  \end{pmatrix}$
  \item $\text{qtr}:V(1) \otimes V(1)^* \to \Q(q,t)^{\tau}$ is the map $v\otimes w^* \to \pi^{p(v)p(w)}q^{-|v|}w^*(v)$. It is given by the matrix
  $\begin{pmatrix}
    q^{-1} & 0 & 0 & \pi q
  \end{pmatrix}$
  \item $\text{coev}:\Q(q,t)^{\tau} \to V(1)^* \otimes V(1) $ given by $1\to qv_{+}\otimes v_{+} + \pi q^{-1} v_{-}\otimes v_{-}$ It is given by the matrix
  $\begin{pmatrix}
     q \\
     0 \\
     0 \\
     \pi q^{-1}
   \end{pmatrix}$
  \item $\text{coqtr}:\Q(q,t)^{\tau} \to V(1) \otimes V(1)^* $ is map $1\to v_{+}\otimes v_{+}+ v_{-}\otimes v_{-}$. It is given by the matrix
  $\begin{pmatrix}
     1 \\
     0 \\
     0 \\
     1
   \end{pmatrix}$
\end{itemize}
as well as the quasi $R$-matrix

$R_{1,1}:V(1)\otimes V(1) \to V(1)\otimes V(1)$ given by the matrix
$\begin{pmatrix}
1 & 0 & 0 & 0 \\
0 & 0 &  q & 0 \\
0 & \pi q & 1-\pi q^2 & 0 \\
0 & 0 & 0 & 1
\end{pmatrix}$

and its inverse $R_{1,1}^{-1}$ given by the matrix
$\begin{pmatrix}
1 & 0 & 0 & 0 \\
0 & 1-\pi q^{-2} & \pi q^{-1} & 0 \\
0 & q^{-1} & 0 & 0 \\
0 & 0 & 0 & 1
\end{pmatrix}$.

So far, we can do everything in this section in the subring $\Q(q)^\pi\subset \Q(q,t)^\tau$, but need to renormalize the maps using $\pi=\tau^2$ to get a link invariant.
In particular, the maps $qtr$ and $R_{1,1}$ pick up a factor of $\tau$, while $coev$ and $R_{1,1}^{-1}$ pick up a factor of $\tau^3$.
The renormalized maps have a diagrammatic interpretation using red strands that can be used to obtain the $U_q({\mf{osp}(1|2)})$ link invariant.

\begin{align}\label{eq:red-maps}
  \hackcenter{ \begin{tikzpicture} [scale=.75]
\draw[thick,   red, ->] (0,1).. controls ++(.25,-.5) and ++(-.25,-.5).. (1,1);
\end{tikzpicture}} & = \tau^3 coev =
\begin{pmatrix}
     \tau^3q \\
     0 \\
     0 \\
     \tau q^{-1}
   \end{pmatrix}  & \qquad \hackcenter{ \begin{tikzpicture} [scale=.75]
\draw[thick,   red, <-] (0,1).. controls ++(.25,-.5) and ++(-.25,-.5).. (1,1);
\end{tikzpicture}} & = coqtr =
\begin{pmatrix}
     1 \\
     0 \\
     0 \\
     1
   \end{pmatrix} \\
\hackcenter{ \begin{tikzpicture} [scale=.75]
\draw[thick,   red, ->] (0,0).. controls ++(.25,.5) and ++(-.25,.5).. (1,0);
\end{tikzpicture}} & = \tau qtr =
\begin{pmatrix}
    \tau q^{-1} & 0 & 0 & \tau^3 q
  \end{pmatrix} & \qquad \hackcenter{ \begin{tikzpicture} [scale=.75]
\draw[thick,   red, <-] (0,0).. controls ++(.25,.5) and ++(-.25,.5).. (1,0);
\end{tikzpicture}} & = ev =
\begin{pmatrix}
    1 & 0 & 0 & 1
  \end{pmatrix} \\
\hackcenter{ \begin{tikzpicture} [scale=.75]
\draw[thick, red, ->] (0,0) to (1,1);
\draw[thick, red, -] (1,0) to (.6,.4);
\draw[thick, red, <-] (0,1) to (.4,.6);
\end{tikzpicture}} &= \tau R_{1,1} =
\begin{pmatrix}
\tau & 0 & 0 & 0 \\
0 & 0 & \tau q & 0 \\
0 & \tau^3q & \tau-\tau^3q^2 & 0 \\
0 & 0 & 0 & \tau
\end{pmatrix}
& \qquad
\hackcenter{ \begin{tikzpicture} [xscale=-1,scale=.75]
\draw[thick,   red, ->] (0,0) to (1,1);
\draw[thick,   red, -] (1,0) to (.6,.4);
\draw[thick,   red, <-] (0,1) to (.4,.6);
\end{tikzpicture}} &= \tau^3 R_{1,1}^{-1} =
\begin{pmatrix}
\tau^3 & 0 & 0 & 0 \\
0 & \tau^3-\tau q^{-2} & \tau q^{-1} & 0 \\
0 & \tau^3q^{-1} & 0 & 0 \\
0 & 0 & 0 & \tau^3
\end{pmatrix}
\end{align}

\begin{theorem}{\cite[Theorem 3.8, Corollary 3.9]{C}}
There is a covariant functor $J$ from the category of oriented tangles to the category of finite dimensional $U_{q,\pi}$-modules that sends a tangle $T$ to the associated $U_{q,\pi}$-module homomorphism defined using \ref{eq:red-maps} and multiplying the resulting homomorphism by $(\pi q)^{wr(T)}$.
In particular, $L$ is an oriented link, then $J_L^1=J(L)$ is Clark's $\mf{osp}(1|2)$-link invariant \cite[Example 3.10]{C}.
\end{theorem}

In \cite[Example 3.10]{C}, Clark showed that the red maps \ref{eq:red-maps} satisfy
\begin{align*}
\hackcenter{ \begin{tikzpicture} [scale=.75]
\draw[thick, red, ->] (0,0) to (1,1);
\draw[thick, red, -] (1,0) to (.6,.4);
\draw[thick, red, <-] (0,1) to (.4,.6);
\end{tikzpicture}}
- q^2\;\;
\hackcenter{ \begin{tikzpicture} [xscale=-1,scale=.75]
\draw[thick, red, ->] (0,0) to (1,1);
\draw[thick, red, -] (1,0) to (.6,.4);
\draw[thick, red, <-] (0,1) to (.4,.6);
\end{tikzpicture}}
\;\;=\;\;
(\tau-\tau^3 q^2)\;\;
\hackcenter{ \begin{tikzpicture} [scale=.75]
\draw[thick,   red,->] (0,0) to (0,1);
\draw[thick,   red,->] (1,0) to (1,1);
\end{tikzpicture}}
\end{align*}
This shows that $J_L^1$ satisfies the skein relation \ref{osp-Skein}.
Then taking the specialization $\tau=t$, which correpsonds to $\pi=-1$, he concludes that the uncolored $U_q(\mf{osp}(1|2))$ link invariant is equal to the $U_{\tau^{-1}q}(\mf{sl}_2)$ link invariant.

\subsection{Representation theoretic definiton of $\hat{J}$}\label{Sec:Our-method}

For the purposes of categorification, we want to define the invariant $\hat{J}$ defined in section 1 using $U_{q,\pi}$-module homomorphisms.
To do this, we need to define $U_{q,\pi}$-module homomorphisms $u: \Q(q,t)^\tau \to V \otimes V$ and $n: V \otimes V \to \Q(q,t)^\tau$ satisfying the relation $\tau R_{1,1}=\tau \Id -q un$.
Then this gives the $\mf{sl}_2$ invariant in the variable $\tau^{-1}q$ as desired.

\begin{definition}
 Let $n:V\otimes V \to \Q(q,t)^{\tau}$ be the odd $\Q(q,t)^\tau$-linear map defined on basis elements by
  \begin{align*}
    v_+\otimes v_+ \to 0 & \qquad v_+ \otimes v_- \to -\tau^3 \\
    v_- \otimes v_- \to 0 & \qquad v_- \otimes v_+ \to \tau^3 q
  \end{align*}
\end{definition}

\begin{theorem}
  $n$ is a odd weight $U_{q,\pi}$-module homomorphism.
\end{theorem}
\begin{proof}
  We need to show that $n(\Delta(E)v\otimes w)=n(\Delta(F)v\otimes w)=0$.
  We know that $\Delta(E)(v_+ \otimes v_+)=\Delta(F)(v_i \otimes v_-)=0$, so it holds in these cases.
  It is also trivial to see that $\Delta(E)(v_{\pm}\otimes v_{\mp}) \in span\{v_+ \otimes v_+\}$, which is sent to $0$ by $n$.
  Likewise, $\Delta(F)(v_{\pm}\otimes v_{\mp})\in span\{v_-\otimes v_-\}$, which is sent to $0$ by $n$.
  Thus, all that remains to check is that $\Delta(F)(v_+ \otimes v_+)$ and $\Delta(E)(v_-\otimes v_-)$ are sent to $0$ by $n$.

  \begin{align*}
    n((E\otimes 1 + JK \otimes E)(v_- \otimes v_-)) &= n(v_+ \otimes v_- +\pi JKv_-\otimes Ev_-)= n(v_+\otimes v_- + q^{-1}v_-\otimes v_+)=0 \\
    n((F\otimes K^{-1} + 1\otimes F)(v_+ \otimes v_+))&=n(q^{-1}v_- \otimes v_+ + v_+ \otimes v_-)=0\\
  \end{align*}
  Therefore, $n(\Delta(x)(v\otimes w))=x\cdot n(v\otimes w) = 0$ for all $x\in \{E,F\}$ and $v,w \in V$, which implies $n$ is a module homomorphism.
\end{proof}

We associate a blue cap with no orientation to this map
\begin{align}\label{n-draw}
  \hackcenter{ \begin{tikzpicture} [scale=.75]
\draw[thick,   blue, -] (0,0).. controls ++(.25,.5) and ++(-.25,.5).. (1,0);
\end{tikzpicture}} & = n =
\begin{pmatrix}
     0 & -\tau^3 & \tau^3q & 0
  \end{pmatrix}
\end{align}

\begin{definition}
  Let $u:\Q(q,t)^\tau \to V \otimes V$ be the odd $\Q(q,t)^\tau$-linear map defined by
  \[
  1\to -\pi q^{-1}v_+ \otimes v_{-} + v_{-}\otimes v_{+}
  \]
\end{definition}

\begin{theorem}
  $u$ is a an odd weight $U_{q,\pi}$-module homomorphism.
\end{theorem}
\begin{proof}
  We need to show that $\Delta(x)u(1)=0$ for $x\in \{E,F\}$.
  \begin{align*}
    \Delta(E)u(1)&=(E\otimes 1)(v_{-}\otimes v_{+})+(JK \otimes E)(-\pi q^{-1}v_+ \otimes v_{-})=v_+\otimes v_+ -\pi q^{-1} JKv_+ \otimes Ev_-\\
    &=v_+\otimes v_+ -v_+\otimes v_+=0
  \end{align*}
  \begin{align*}
    \Delta(F)u(1)&=(F\otimes K^{-1})(-\pi q^{-1}v_+ \otimes v_{-})+(1\otimes F)(v_-\otimes v_+)=-\pi q^{-1}v_-\otimes K^{-1}v_-+\pi v_-\otimes v_- \\
    &=-\pi v_-\otimes v_- +\pi v_- \otimes v_-=0
  \end{align*}
\end{proof}

We associate a cup without orientation to the map $u$

\begin{align}\label{blue-cup-draw}
  \hackcenter{ \begin{tikzpicture} [scale=.75]
\draw[thick,   blue, -] (0,1).. controls ++(.25,-.5) and ++(-.25,-.5).. (1,1);
\end{tikzpicture}} & = u =
\begin{pmatrix}
     0\\
     -\pi q^{-1}\\
     1 \\
     0
  \end{pmatrix}
\end{align}

Then $nu=\tau q^{-1}+\tau^3 q$.

Define $E'= un: V\otimes V\to V\otimes V$. $E'$ is the even map
\begin{align*}
  \hackcenter{ \begin{tikzpicture} [scale=.75]
\draw[thick,   blue] (0,0).. controls ++(.25,.5) and ++(-.25,.5).. (1,0);
\draw[thick,   blue] (0,1).. controls ++(.25,-.5) and ++(-.25,-.5).. (1,1);
\end{tikzpicture}} &= E'= un =
\begin{pmatrix}
     0\\
     -\pi q^{-1}\\
     1 \\
     0
  \end{pmatrix}
\begin{pmatrix}
     0 & -\tau^3 & \tau^3q & 0
  \end{pmatrix}
=
\begin{pmatrix}
  0 & 0 & 0 & 0 \\
  0 & \tau q^{-1} & -\tau & 0 \\
  0 & -\tau^3 & \tau^3 q  & 0 \\
  0 & 0 & 0 & 0
\end{pmatrix}
\end{align*}

We can derive all the relations in \eqref{bracket-poly2-1}, \eqref{bracket-poly2-2}, and \eqref{bracket-poly2-3} with this understanding of the cap and cup as $U_{q,\pi}$-module homomorphisms.

\begin{align*}
\hackcenter{ \begin{tikzpicture} [scale=.75]
\draw[thick, red, -] (0,0) to (1,1);
\draw[thick, red, -] (1,0) to (.6,.4);
\draw[thick, red, -] (0,1) to (.4,.6);
\end{tikzpicture}} \;=
 \begin{pmatrix}
\tau & 0 & 0 & 0 \\
0 & 0 & \tau q & 0 \\
0 & \tau^3q & \tau-\tau^3q^2 & 0 \\
0 & 0 & 0 & \tau
\end{pmatrix}&=
\tau Id-q
\begin{pmatrix}
  0 & 0 & 0 & 0 \\
  0 & \tau q^{-1} & -\tau & 0 \\
  0 & -\tau^3 & \tau^3 q  & 0 \\
  0 & 0 & 0 & 0
\end{pmatrix}=
\tau\;\;
\hackcenter{ \begin{tikzpicture} [scale=.75]
\draw[thick,   blue] (0,0) to (0,1);
\draw[thick,   blue] (1,0) to (1,1);
\end{tikzpicture}}
-q\;\;
\hackcenter{ \begin{tikzpicture} [scale=.75]
\draw[thick,   blue] (0,0).. controls ++(.25,.5) and ++(-.25,.5).. (1,0);
\draw[thick,   blue] (0,1).. controls ++(.25,-.5) and ++(-.25,-.5).. (1,1);
\end{tikzpicture}}
\end{align*}

\begin{align*}
\hackcenter{ \begin{tikzpicture} [xscale=-1,scale=.75]
\draw[thick, red, -] (0,0) to (1,1);
\draw[thick, red, -] (1,0) to (.6,.4);
\draw[thick, red, -] (0,1) to (.4,.6);
\end{tikzpicture}} \;=
  \begin{pmatrix}
\tau^3 & 0 & 0 & 0 \\
0 & \tau^3-\tau q^{-2} & \tau q^{-1} & 0 \\
0 & \tau^3q^{-1} & 0 & 0 \\
0 & 0 & 0 & \tau^3
\end{pmatrix}&=
\tau^3 Id- q^{-1}
\begin{pmatrix}
  0 & 0 & 0 & 0 \\
  0 & \tau q^{-1} & -\tau & 0 \\
  0 & -\tau^3 & \tau^3 q  & 0 \\
  0 & 0 & 0 & 0
\end{pmatrix}
=
\tau^3\;\;
\hackcenter{ \begin{tikzpicture} [scale=.75]
\draw[thick,   blue] (0,0) to (0,1);
\draw[thick,   blue] (1,0) to (1,1);
\end{tikzpicture}}
-q^{-1}\;\;
\hackcenter{ \begin{tikzpicture} [scale=.75]
\draw[thick,   blue] (0,0).. controls ++(.25,.5) and ++(-.25,.5).. (1,0);
\draw[thick,   blue] (0,1).. controls ++(.25,-.5) and ++(-.25,-.5).. (1,1);
\end{tikzpicture}}
\end{align*}
\begin{align*}
   \hackcenter{ \begin{tikzpicture} [scale=.75]
\draw[thick,   blue] (0,0).. controls ++(.25,.5) and ++(-.25,.5).. (1,0);
\draw[thick,   blue] (0,0).. controls ++(.25,-.5) and ++(-.25,-.5).. (1,0);
\end{tikzpicture}} &=nu=
\begin{pmatrix}
     0 & -\tau^3 & \tau^3q & 0
  \end{pmatrix}
\begin{pmatrix}
     0\\
     -\pi q^{-1}\\
     1 \\
     0
  \end{pmatrix}=\tau q^{-1}+\tau^3 q
\end{align*}

\begin{align*}
    \hackcenter{ \begin{tikzpicture} [xscale=-1, scale=.75]
%
\draw[thick,   blue] (0,0).. controls ++(.25,.5) and ++(-.25,.5).. (1,0);
\draw[thick,   blue] (1,0).. controls ++(.25,-.5) and ++(-.25,-.5).. (2,0);
\draw[thick, blue] (2,0) to (2,1);
\draw[thick, blue] (0,0) to (0,-1);
\end{tikzpicture}}&=(1\otimes n)\circ (u\otimes 1):V\to V
\end{align*}
This is the map sending
\begin{align*}
  v_{\pm} &\overset{u\otimes 1}{\rightarrow}-\pi q^{-1}v_+\otimes v_- \otimes v_{\pm}+ v_-\otimes v_+\otimes v_{\pm} \overset{1\otimes n}{\rightarrow} -\pi q^{-1}v_+\otimes n(v_- \otimes v_{\pm})+\pi  v_-\otimes n(v_+\otimes v_{\pm})=-\tau v_{\pm}.
\end{align*}
Hence,
\begin{align*}
    \hackcenter{ \begin{tikzpicture} [xscale=-1, scale=.75]
%
\draw[thick,   blue] (0,0).. controls ++(.25,.5) and ++(-.25,.5).. (1,0);
\draw[thick,   blue] (1,0).. controls ++(.25,-.5) and ++(-.25,-.5).. (2,0);
\draw[thick, blue] (2,0) to (2,1);
\draw[thick, blue] (0,0) to (0,-1);
\end{tikzpicture}}&=-\tau \;\;
\hackcenter{ \begin{tikzpicture} [xscale=-1, scale=.75]
\draw[thick, blue] (1,1) to (1,-1);
\end{tikzpicture}}
\end{align*}

The other zig-zag straightens up to a different power of $\tau$.
\begin{align*}
   \hackcenter{ \begin{tikzpicture} [ scale=.75]
%
\draw[thick,   blue] (0,0).. controls ++(.25,.5) and ++(-.25,.5).. (1,0);
\draw[thick,   blue] (1,0).. controls ++(.25,-.5) and ++(-.25,-.5).. (2,0);
\draw[thick, blue] (2,0) to (2,1);
\draw[thick, blue] (0,0) to (0,-1);
\end{tikzpicture}} & =(n\otimes 1)\circ (1\otimes u):V \to V
\end{align*}
Note that $(1\otimes u):V \to V\otimes V\otimes V$ is the map sending
\[v_{\pm} =v_{\pm}\otimes 1 \to \pi^{p(v_{\pm})}v_{\pm}\otimes u(1)=\pi^{p(v_{\pm})}(-\pi q^{-1}  v_{\pm}\otimes  v_+\otimes v_- + v_{\pm}\otimes  v_-\otimes v_+).\]
Thus, $(n\otimes 1)\circ (1\otimes u)$ is the map sending
\begin{align*}
  v_+ &\overset{1\otimes u}{\rightarrow} (-\pi q^{-1}  v_{+}\otimes  v_+\otimes v_- + v_{+}\otimes  v_-\otimes v_+)\overset{n\otimes 1}{\rightarrow} -\tau^3 v_+  \\
  v_- &\overset{1\otimes u}{\rightarrow}  \pi(-\pi q^{-1}  v_{-}\otimes  v_+\otimes v_- + v_{-}\otimes  v_-\otimes v_+)\overset{n\otimes 1}{\rightarrow} -\tau^3 v_-.
\end{align*}
Hence,
\begin{align*}
    \hackcenter{ \begin{tikzpicture} [xscale=1, scale=.75]
%
\draw[thick,   blue] (0,0).. controls ++(.25,.5) and ++(-.25,.5).. (1,0);
\draw[thick,   blue] (1,0).. controls ++(.25,-.5) and ++(-.25,-.5).. (2,0);
\draw[thick, blue] (2,0) to (2,1);
\draw[thick, blue] (0,0) to (0,-1);
\end{tikzpicture}}&=-\tau^3 \;\;
\hackcenter{ \begin{tikzpicture} [xscale=1, scale=.75]
\draw[thick, blue] (1,1) to (1,-1);
\end{tikzpicture}}
\end{align*}

Therefore, we also have the relation
\begin{align*}
    \hackcenter{ \begin{tikzpicture} [scale=.75]
\draw[thick,   blue] (0,1).. controls ++(.25,-.5) and ++(-.25,-.5).. (1,1);
\draw[thick,   blue] (1,1).. controls ++(.25,.5) and ++(-.25,.5).. (2,1);
\draw[thick,   blue] (0,0).. controls ++(.25,.5) and ++(-.25,.5).. (1,0);
\draw[thick,   blue] (1,0).. controls ++(.25,-.5) and ++(-.25,-.5).. (2,0);
\draw[thick, blue] (2,0) to (2,1);
\draw[thick, blue] (0,0) to (0,-1);
\draw[thick, blue] (0,1) to (0,2);
\end{tikzpicture}}
\;\;=\;\;
  \hackcenter{ \begin{tikzpicture} [xscale=-1, scale=.75]
\draw[thick,   blue] (0,1).. controls ++(.25,-.5) and ++(-.25,-.5).. (1,1);
\draw[thick,   blue] (1,1).. controls ++(.25,.5) and ++(-.25,.5).. (2,1);
\draw[thick,   blue] (0,0).. controls ++(.25,.5) and ++(-.25,.5).. (1,0);
\draw[thick,   blue] (1,0).. controls ++(.25,-.5) and ++(-.25,-.5).. (2,0);
\draw[thick, blue] (2,0) to (2,1);
\draw[thick, blue] (0,0) to (0,-1);
\draw[thick, blue] (0,1) to (0,2);
\end{tikzpicture}}
\;\;=\;\;
\hackcenter{ \begin{tikzpicture} [scale=.75]
\draw[thick,   blue] (0,-1) to (0,2);
\end{tikzpicture}}
\end{align*}

Therefore, $\hat{J}$ can be understood using representation theory of $U_{q,\pi}$.

\begin{definition}
  Let $d:=\tau^3 q+\tau q^{-1}$. Define $\sigma_i\in B_n$ using the positive crossing involving strands $i$ and $i+1$.
  Then define $\rho_{q,\tau}:B_n \to TL_n(d)$ given by
  \[
  \rho_{q,\tau}(\sigma_i)=(\tau Id -q E_i).
  \]
\end{definition}
\begin{theorem}
  $\rho_{q,\tau}:B_n \to TL_n(d)$ defines a $B_n$-representation inside $TL_n(d)$, namely a group homomorphism from $B_n$ to the group of invertible elements of $TL_n(d)$.
\end{theorem}

\section{Two Categorifications}\label{section:cat}
As we know from \cite{C,Blu09}, the $\mf{osp}(1|2)$ polynomial (and hence our invariant $\hat{J}$) is practically identical to the Jones polynomial.
Now that we have this Kauffman-bracket type definition of $\hat{J}$, it should not come as much of a surprise that one can categorify this invariant using a homology theory that is practically identical to Khovanov homology.
Indeed, we show in this section that one can repeat the construction of Khovanov homology seen in \cite{BN2} but where an extra $\Z_4$-grading is attached to everything to get a categorification.
We then construct a $\Z\times \Z_4$-graded version of Putyra's covering Khovanov homology \cite{P} in a similar fashion.

\begin{notation}
  An "e-graded" object refers to an object with a $\Z\times \Z_4$-grading.
\end{notation}
For example, an e-graded vector space is a $\Z\times \Z_4$-graded vector space.
\begin{definition}
  Let $\cat{egVec}$ be subcategory of $\cat{Vect}$ consisting of $\Z\times \Z_4$-graded vector spaces and degree-preserving linear maps.
\end{definition}
This is a symmetric monoidal category with the same braiding as $\cat{Vect}$.
The tensor product of linear maps is given by $(f\otimes g)(v\otimes w)=f(v)\otimes g(w)$.
We bring this example up to emphasize that the $\Z_4$-grading doesn't impact the tensor product of linear maps or the braiding like the $\Z_2$-grading does for superspaces.
For the remainder of the paper, we choose to work with modules instead of vector spaces.

An e-graded module $M$ over a e-graded ring $R$ is a $\Z\times \Z_4$-graded $R$-module such that $R_{(a,b)}M_{(x,y)}\subset M_{(x+a,y+b)}$.
We refer to the degree of a homogeneous element in some e-graded module with $deg(-)$ for the rest of the paper.
\begin{definition}
  The category of e-graded $R$-modules $\cat{R-emod}$ is the subcategory of $\cat{R-mod}$ whose objects are e-graded $R$-modules and whose morphisms are degree preserving $R$-module homomorphisms.
\end{definition}
For our discussion on Khovanov homology, we set $R=\Z$ (concentrated in degree $(0,0)$).
\begin{definition}
  Let $W=\bigoplus\limits_{(m,n)\in \Z \times \Z_4} W_{m,n}$ be a e-graded $R$-module with homogeneous components $W_{m,n}$. The graded dimension of $W$ is the power series $\text{dim}_{q,\tau}(W)=\sum_{m,n}q^m\tau^n \text{dim}(W_{m,n})$. When $W$ is finite dimensional, this lies in $\Z^\tau[q,q^{-1}]$.
\end{definition}
\begin{definition}
  Let $W= \bigoplus W_{m,n}$ as before. Then define $\{(a,b)\}$ to be the degree shift. That is, we set $W\{(a,b)\}_{(m,n)}=W_{(m-a,n-b)}$ so that $\text{dim}_{q,\tau}W\{(a,b)\}=q^{a}\tau^b \text{dim}_{q,\tau} W$.
\end{definition}

Let $\cat{2Cob}$ denote the category of 2D-cobordisms $Cob^3(\emptyset)$ from \cite{BN2}.
The objects are closed 1-dimensional manifolds (collections of circles), and the morphisms are given by 2-dimensional cobordisms between these collections of circles up to boundary preserving isotopies.
We also use the grading conventions in \cite[Section 6]{BN2}, which tell us that the degree of a cobordism $W: \Sigma_1\{m\}\to \Sigma_2\{n\}$ is given by $\chi(W)+n-m$.
We will later attach an additional $\Z_4$-grading to this category.


\subsection{Normal Khovanov Homology}
We briefly recall how the complex $Kh(L)$ for Khovanov homology of an oriented link $L$ with is defined.
We use notation similar to \cite{BN2, P} to define the Khovanov complexes so that we can easily reference it again when we look at a similar e-graded version of covering Khovanov homology.
Let $X$ be a crossing in $L$, and $L_0, L_1$ denote the diagram obtained by replacing $X$ with its $0$ and $1$-resolution respectively.
Replace each crossing $X$ with a two term chain complex $D_X: L_0 \to L_1\{1\}$, where $L_0$ lives in cohomological degree $0$ is $X$ is a positive crossing and lives in cohomological degree $-1$ otherwise.
$D_X$ is a saddle cobordism, so the degree of $D_X$ is $1-1=0$ (see \cite[Section 6]{BN2}).

Doing this for all the crossings gives the cube of resolutions with a saddle cobordism associated to each edge.
Then for any sign assignment $\epsilon$ for the edges of the cube of resolutions so that faces anti-commute (e.g \cite[Section 3.2]{BN1}), one obtains a chain complex of cobordisms $[L]'$ defined by
\begin{align*}
  [L]_{r-n_-}'&= \bigoplus_{\sum z_i =r} L_z'\{r\} \\
  d|_{L_z'}&=\sum\limits_{\alpha: z\to w}\epsilon(\alpha)D_{\alpha}
\end{align*}
where $\alpha$ sums over edges out of $z$.
Let $j$ be the entry where $z$ and $w$ differ and $X_j$ denote the $j$-th crossing; then $D_\alpha= D_{X_j}$.
Then $Kh(L)^i=[L]_{i}'\{n_+-2n_-\}$ is the Khovanov complex.
One then applies the usual TQFT to get a complex of $\Z$-modules whose homology is Khovanov homology.

\subsection{e-graded Khovanov Homology}\label{section: e-Khov}
We construct the e-graded version of Khovanov homology in a nearly identical manner.
We can define an e-grading on $\cat{2Cob}$ by adding formal objects $\Sigma\{(a,b)\}$ for every $1$-manifold $\Sigma$ and $(a,b)\in \Z \times \Z_4$, and defining the degree $deg(W)$ of a cobordism $W:\Sigma_0\{(m,m')\}\to \Sigma_{1}\{(n,n')\}$ by
\[
deg(W)=(\chi(W)+n-m, -\chi(W)+n'-m')
\]
Pants cobordisms have degree $(-1,1)$ and caps and cups have degree $(1,-1)$.

Given an oriented link diagram $L$, construct a chain complex $[L]$ of cobordisms as follows.
Replace each crossing of $L$ with the chain complex
\begin{align}\label{crossing-complexes}
\left[
  \hackcenter{ \begin{tikzpicture} [scale=.75]
\draw[thick,   blue] (0,0) to (1,1);
\draw[thick,   blue] (1,0) to (.6,.4);
\draw[thick,   blue] (0,1) to (.4,.6);
\end{tikzpicture}}\right]
&=
\underline{
\hackcenter{ \begin{tikzpicture} [scale=.75]
\draw[thick,   blue] (0,0) to (0,1);
\draw[thick,   blue] (1,0) to (1,1);
\end{tikzpicture}}} \{(0,1)\}
\to
\hackcenter{ \begin{tikzpicture} [scale=.75]
\draw[thick,   blue] (0,0).. controls ++(.25,.5) and ++(-.25,.5).. (1,0);
\draw[thick,   blue] (0,1).. controls ++(.25,-.5) and ++(-.25,-.5).. (1,1);
\end{tikzpicture}}
\{(1,0)\}
 \\
 \left[
 \hackcenter{ \begin{tikzpicture} [xscale=-1, scale=.75]
\draw[thick,   blue] (0,0) to (1,1);
\draw[thick,   blue] (1,0) to (.6,.4);
\draw[thick,   blue] (0,1) to (.4,.6);
\end{tikzpicture}}\right]
&=
\hackcenter{ \begin{tikzpicture} [scale=.75]
\draw[thick,   blue] (0,0).. controls ++(.25,.5) and ++(-.25,.5).. (1,0);
\draw[thick,   blue] (0,1).. controls ++(.25,-.5) and ++(-.25,-.5).. (1,1);
\end{tikzpicture}}
\{(-1,0)\}
\to
\underline{
\hackcenter{ \begin{tikzpicture} [scale=.75]
\draw[thick,   blue] (0,0) to (0,1);
\draw[thick,   blue] (1,0) to (1,1);
\end{tikzpicture}}} \{(0,-1)\}
\end{align}
where the underlined term in each complex lives in cohomological degree $0$.
The maps between the $0$ and $1$-resolution in each complex are given by a saddle cobordism just like before.

Then the complex $[L]$ is defined by
\begin{align}\label{pre-e-Khov}
  [L]^{r-n_-} &= \bigoplus_{\sum z_i=r} L_z' \{(-n_-+r, n_+-r)\} \\
  d|_{L_z'} &= \sum\limits_{\alpha: z \to w}\epsilon(\alpha)D_{\alpha}
\end{align}
where $D_\alpha=D_{X_j}$ is one of the two maps shown above.
We then get the e-graded version of the Khovanov complex by shifting the degree of each term.
\begin{definition}
The chain complex $eKh(L)$ is obtained by shifting the degree of each chain group by $\{wr(L), 2wr(L)\}$
  \begin{align}\label{e-Khov-1}
  eKh(L)&=[L]\{(wr(L), 2wr(L)\}
\end{align}
\end{definition}

Let $A=R[X]/(X^2)$ be the free e-graded $R$-module of rank $2$ spanned by $\1$ and $X$ in degrees $deg(\1)=(1,-1)$ and $deg(X)=(-1,1)$.
We equip $A$ with the following maps to make it a Frobenius algebra.
\begin{align}\label{Frob-str-1}
  \varepsilon: A\to R & \left\{\begin{matrix}
                                \1 \to 0 \\
                                X\to 1
                              \end{matrix} \right.\\
  \eta: R \to A & \left\{ 1\to \1 \right. \\
  m:A\otimes A \to A & \left\{  \begin{matrix}
                                  \1 \otimes \1 \to \1 & \1\otimes X \to X \\
                                  X\otimes \1 \to X & X\otimes X \to 0
                                \end{matrix}  \right. \\
  \Delta: A\to A\otimes A & \left\{ \begin{matrix}
                                      \1 \to \1\otimes X + X \otimes \1 \\
                                      X \to X \otimes X
                                    \end{matrix}\right.
\end{align}
We have $deg(m)=deg(\Delta)=(-1,1)$ and $deg(\varepsilon)=deg(\eta)=(1,-1)$.

\begin{definition}
  Let $\mathcal{F}:\cat{2Cob} \to \cat{R-emod}$ be the TQFT corresponding to the Frobenius algebra $A$. Denote $M_n$ as the object in $\cat{2Cob}$ consisting of $n$-disjoint circles. $\mathcal{F}$ sends a collection of $k$ circles to $A^{\otimes k}$ and
  \begin{itemize}
    \item The two pair-of-pants cobordisms $M_2\to M_1$ and $M_1\to M_2$ to $m$ and $\Delta$ respectively.
    \item The birth $M_0 \to M_1$ and death $M_1 \to M_0$ to $\eta$ and $\varepsilon$ respectively.
  \end{itemize}
\end{definition}
It is easy to see that $\mathcal{F}$ is well defined and degree preserving.
Apply the TQFT $\mathcal{F}:\cat{2Cob} \to \cat{R-emod}$ to turn the chain complex $eKh[L]$ into a chain complex $C(L):=\mathcal{F}(eKh[L])$ of e-graded $R$-modules.

\begin{theorem}\label{homology-link-invariant}
  The homology $H(C(L))$ of the chain complex $C(L)$ is an invariant of oriented links.
\end{theorem}

To show this, we need to show that the Reidemeister moves hold.
Since forgetting the $\Z_4$-grading recovers Khovanov homology exactly, the arguments from \cite{BN1} for why Khovanov homology is a link invariant can be used to prove this theorem.
To demonstrate, we prove that this homology satisfies R1 as a testament to how the arguments for Khovanov homology work verbatim for this homology theory.

\subsubsection{R1}
\[
C:=C
\left(
\hackcenter{ \begin{tikzpicture} [scale=.75]
\draw[thick,   blue] (0,0) to (1,1);
\draw[thick,   blue] (1,0) to (.6,.4);
\draw[thick,   blue] (0,1) to (.4,.6);
\draw[thick, blue] (1,1) .. controls ++(.5,.25) and ++(.5, -.25) ..(1,0);
\end{tikzpicture}}
\right)=\left[0\to
 \left(A\otimes C\left(
\hackcenter{ \begin{tikzpicture} [scale=.75]
\draw[thick,   blue] (0,0) to (0,1);
\end{tikzpicture}}\right)
\right)\{(0,1)\} \overset{m}{\to}
C\left(
\hackcenter{ \begin{tikzpicture} [scale=.75]
\draw[thick,   blue] (0,0) to (0,1);
\end{tikzpicture}}\right)
\{(1,0)\}
\to 0 \right]\{(1,2)\}
\]
where the first non-zero term is in homological degree 0.

Recall that $m(\1 \otimes a)=a$, and that $A=R\1 \oplus R X$, so $C$ has an acyclic subcomplex $C'$.
\[
C'=
\left[0\to
\left(R \1 \otimes C\left(
\hackcenter{ \begin{tikzpicture} [scale=.75]
\draw[thick,   blue] (0,0) to (0,1);
\end{tikzpicture}}\right)
\right)\{(0,1)\} \overset{m}{\to}
C\left(
\hackcenter{ \begin{tikzpicture} [scale=.75]
\draw[thick,   blue] (0,0) to (0,1);
\end{tikzpicture}}\right)
\{(1,0)\}
\to 0
\right]\{(1,2)\}
\]

Taking the quotient of $C$ by this subcomplex, we get
\[
C/C'\cong
\left[
0\to \left(R X \otimes C\left(
\hackcenter{ \begin{tikzpicture} [scale=.75]
\draw[thick,   blue] (0,0) to (0,1);
\end{tikzpicture}}\right)
\right)\{(0,1)\}\to
0
\to 0
\right]\{(1,2)\}
\cong
C\left(
\hackcenter{ \begin{tikzpicture} [scale=.75]
\draw[thick,   blue] (0,0) to (0,1);
\end{tikzpicture}}\right)
\]
So $H(C)=H(C/C')$.
The other R1 relation holds by a similar argument.

\subsubsection{Categorification of $\hat{J}$}
\begin{definition}
  Let $C$ be a chain complex and denote the term in cohomological degree $r$ by $C^r$. Then define the e-graded Euler characteristic $\chi_{q,\tau}$ of this complex by
  \begin{align*}
  \chi_{q,\tau}(C)=\sum\limits_{r\in \Z} (-1)^{r} \text{dim}_{q,\tau}(C^r).
\end{align*}
\end{definition}

\begin{theorem}
  $\hat{J}(L)= \chi_{q,\tau}(C(L))$
\end{theorem}

\begin{proof}
   If $D$ is a collection of $k$ disjoint circles, then $C(D)=0\to A^{\otimes k}\to 0$. So $\chi_{q,\tau}(C(D))=(\tau^{-1}q+\tau q^{-1})^k$.

   Let $D_1,D_2, D_3$ be diagrams that differ as shown below
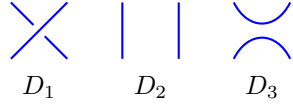
\[
\hackcenter{ \begin{tikzpicture} [scale=.75]
\draw[thick,   blue] (0,0) to (1,1);
\draw[thick,   blue] (1,0) to (.6,.4);
\draw[thick,   blue] (0,1) to (.4,.6);
\node at (.5,-.5) {$D_1$};
\end{tikzpicture}}
\qquad
\hackcenter{ \begin{tikzpicture} [scale=.75]
\draw[thick,   blue] (0,0) to (0,1);
\draw[thick,   blue] (1,0) to (1,1);
\node at (.5,-.5) {$D_2$};
\end{tikzpicture}}
\qquad
\hackcenter{ \begin{tikzpicture} [scale=.75]
\draw[thick,   blue] (0,0).. controls ++(.25,.5) and ++(-.25,.5).. (1,0);
\draw[thick,   blue] (0,1).. controls ++(.25,-.5) and ++(-.25,-.5).. (1,1);
\node at (.5,-.5) {$D_3$};
\end{tikzpicture}}
\]

Then $C(D_1)$ is isomorphic, up to a shift, to the cone of a map of complexes $C(D_2)\{(1,3)\}\to C(D_3)\{(2,2)\}$.
Therefore, \[\chi(C(D_1))=\pi q ( \tau \chi_{q,\tau}(C(D_2))- q \chi_{q,\tau}(C(D_3))),\] which matches with $\hat{J}$.
\end{proof}
This shows that $H(C(L))$ is a categorification of the link invariant $\hat{J}$.

\subsection{e-graded Covering Khovanov Homology}
Let $R=\Z^\pi$ for the remainder of the paper.
This section assumes prior knowledge of Putyra's work on Covering Khovanov homology \cite{P}.
For the sake of brevity we omit most of the technical details and definitions involving "chronology".
Covering Khovanov homology constructs a chain complex of modules over the ring $\mathbf{k}:=\Z[X,Y,Z]/(X^2=Y^2=1)$ such that the homology of the chain complex if we specialize $X,Y,Z$ to be 1 is regular Khovanov homology and the homology of the chain complex where we take $X, Z$ to be $1$ but $Y$ to be $-1$ is Odd Khovanov homology.
For our purposes, we can describe this as a chain complex of $R=\Z^\pi$ modules with $\pi$ playing the role of $Y$.

We first construct an e-graded version of the generalized Khovanov bracket $eKh_{cov}(L)$ from \cite[Definition 5.9]{P} by following the exact same procedure as Putyra except we handle the gradings and shifts exactly as in Section \ref{section: e-Khov}.
This is a chain complex in the additive closure of the $R$-linear category of chronological cobordisms $\cat{RChCob}$.

Then one can show that this complex is an invariant of oriented links if we consider it as a chain complex of chronological cobordisms modulo the S, T, and 4Tu relations \cite[Section 7]{P}.
The quotient of $\cat{RChCob}$ by these relations is denoted $\cat{RChCob}_{/\ell}$.
The proof of invariance is exactly the same as the proof presented in \cite[Section 7]{P}.

\begin{definition}{\cite[Definition 10.1]{P}}
  Set $G=\Z\times \Z_4 \times \Z_2$, and $U(R)$ to be the group of invertible elements in $R$. Then define a group homomorphism $\lambda: G \times G \to U(R)$ by
  \[
  \lambda((a,b,c)(d,e,f))=\pi^{cf}
  \]
  We also define a e-graded super tensor product of $G$-graded modules as usual, but for homogeneous homomorphisms $f,g$ we define $f\otimes g$ by the formula:
  \[
  (f\otimes g)(m\otimes n) = \lambda(deg(g), deg(m))f(m)\otimes g(n)
  \]
  This defines a symmetric monoidal structure on the category of $G$-graded $R$-modules, which we call $\cat{R-egSmod}$.
\end{definition}
We refer to $\cat{R-egSmod}$ as the category of e-graded $R$-supermodules.
\begin{definition}
  Let $A=R[X]/(X^2)$ be the free e-graded $R$-supermodule of rank 2 spanned by $\1$ and $X$ in degrees $deg(\1)=(1,-1,0)$ and $deg(X)=(-1,1,1)$.
  Equip it with the following operations
  \begin{align}\label{Frob-str-1}
  \varepsilon: A\to R & \left\{\begin{matrix}
                                \1 \to 0 \\
                                X\to 1
                              \end{matrix} \right.\\
  \eta: R\to A & \left\{ 1\to \1 \right. \\
  m:A\otimes A \to A & \left\{  \begin{matrix}
                                  \1 \otimes \1 \to \1 & \1\otimes X \to X \\
                                  X\otimes \1 \to X & X\otimes X \to 0
                                \end{matrix}  \right. \\
  \Delta: A\to A\otimes A & \left\{ \begin{matrix}
                                      \1 \to \pi \1\otimes X + X \otimes \1 \\
                                      X \to X \otimes X
                                    \end{matrix}\right.
\end{align}
These maps have degrees
\[deg(\Delta)=(-1,1,1), \; deg(m)=(-1,1,0), \;deg(\varepsilon)=(1,-1,1), \; deg(\eta)=(1,-1,0).\]
\end{definition}

The pair $(R,A)$ with these operations defines a Chronological Frobenius system as in \cite[Definition 10.5]{P}
Then we obtain a well defined, degree preserving chronological TQFT functor $\mathcal{F}:\cat{RChCob}\to \cat{R-egSmod}$ as in \cite[Proposition 10.6]{P}.
This functor preserves the S, T, and 4Tu relations, so we get a well-defined functor $\mathcal{F}_{cov}:\cat{RChCob}_{/\ell} \to \cat{R-egSmod}$.

Let $C_{cov}(L)$ be the chain complex of e-graded $R$-supermodules given by $C_{cov}(L):=\mathcal{F}_{cov}(eKh_{cov}(L))$.
When $\pi=1$, $C_{cov}(L)=C(L)$ by construction, which is the chain complex for e-graded Khovanov homology.
We call the homology $H(C_{cov}(L))$ of this complex the \emph{e-graded version of covering Khovanov homology}.
If we forget the $\Z_4$-grading of each module in the complex, then we obtain the covering Khovanov homology $\mathcal{H}_{cov}(L)$ defined in \cite[Example 10.7]{P}.

\bibliographystyle{plain}
\bibliography{bib_osp_1}

\end{document}